\documentclass[letterpaper, 10 pt, conference]{ieeeconf}  

\IEEEoverridecommandlockouts 

\usepackage{afterpage}
\usepackage{epsfig} 
\usepackage{times} 
\usepackage{amsmath} 
\usepackage{amssymb}  
\usepackage{amsfonts}
\usepackage{mathptmx} 
\usepackage{tabularx}
\usepackage{epstopdf}
\usepackage{subfigure}
\usepackage{float}
\usepackage{epic,color}
\usepackage{mathrsfs}
\usepackage{dsfont,MnSymbol}
\usepackage{algorithm}
\usepackage{algorithmic}
\usepackage{multirow} 
\usepackage[T1]{fontenc}

\newcounter{rmnum}

\newcounter{anum}

\def\IEEEQEDclosed{\mbox{\rule[0pt]{1.3ex}{1.3ex}}}

\def\qed{\ifmmode\IEEEQEDclosed\else{\unskip\nobreak\hfil
\penalty50\hskip1em\null\nobreak\hfil\IEEEQEDclosed
\parfillskip=0pt\finalhyphendemerits=0\endgraf}\fi}

\def\qed{\hspace*{\fill}~\IEEEQED\par\endtrivlist\unskip}

\def\Re{\mathbb{R}}

\def\Sec#1{Sec.~\ref{#1}}

\def\Appendix#1{Appendix.~\ref{#1}}
\def\notes#1{\marginpar{\tiny #1}\typeout{Notes!
Notes!
Notes!
}}
\renewcommand{\notes}[1]{\typeout{notes!}}
\def\FRAC#1#2#3{\genfrac{}{}{}{#1}{#2}{#3}}
\def\half{{\mathchoice{\FRAC{1}{1}{2}}%
{\FRAC{2}{1}{2}}%
{\FRAC{3}{1}{2}}%
{\FRAC{4}{1}{2}}}}

\newcommand{\tr}{\mbox{tr}}

\def\Re{\field{R}}

\def\Sec#1{Sec.~\ref{#1}}
\def\eqdef{\mathrel{:=}}

\def\transpose{{\hbox{\rm\tiny T}}}
\def\clB{{\cal B}}

\def\clZ{{\cal Z}}

\def\Sec#1{Sec~\ref{#1}}

\def\E{{\sf E}}

\def\spm#1{\notes{\textcolor{blue}{SPM: #1}}}

\def\Sec#1{Sec.~\ref{#1}}

\def\IEEEQEDclosed{\mbox{\rule[0pt]{1.3ex}{1.3ex}}}
\def\qed{\nobreak\hfill\IEEEQEDclosed}

\newcommand{\innov}{{\sf I}}

\def\clZ{{\cal Z}}

\def\eqdef{\mathrel{:=}}

\newtheorem{theorem}{Theorem}

\newtheorem{remark}{Remark}
\newtheorem{proposition}{Proposition}

\def\beq{\begin{eqnarray}} 
\def\bc{\begin{center}} 
\def\be{\begin{enumerate}}
\def\bi{\begin{itemize}} 
\def\bs{\begin{small}}
\def\bS{\begin{slide}}
\def\ec{\end{center}} 
\def\ee{\end{enumerate}}
\def\ei{\end{itemize}}
\def\es{\end{small}}
\def\eS{\end{slide}}
\def\eeq{\end{eqnarray}}

\newcommand{\newP}[1]{\medskip\noindent{\bf #1:}}

\newcommand{\ud}{\,\mathrm{d}}

\def\Re{\mathbb{R}}
\def\E{{\sf E}}

\def\Sec#1{Sec.~\ref{#1}}

\def\Prop#1{Prop.~\ref{#1}}

\def\clB{{\cal B}}

\def\clZ{{\cal Z}}

\renewcommand{\Re}{\mathbb{R}}

\def\eqdef{\mathbin{:=}}

\def\FRAC#1#2#3{\genfrac{}{}{}{#1}{#2}{#3}}

\def\bS{\mathbb{S}}

\def\jin#1{\notes{\textcolor{red}{JIN: #1}}}

\def\dv{\operatorname{diag}}

\def\clU{{\cal U}}
\def\clA{{\cal A}}

\def\sJ{{\sf J}}

\def\lsq#1#2{L_{#1}^2([0,T]\,;{#2})}

\newlength{\noteWidth}
\long\def\notes#1{\ifinner
             {\tiny #1}
             \else
              \marginpar{\parbox[t]{\noteWidth}{\raggedright\tiny #1}}
               \fi}

\long\def\notes#1{}

\title{\LARGE \bf What is the Lagrangian for Nonlinear Filtering?}

\author{Jin-Won Kim, Prashant G. Mehta and Sean P. Meyn
\thanks{Financial support from the 
NSF grant 1761622 and the 
ARO grant W911NF1810334 is gratefully acknowledged. 
}
\thanks{J-W. Kim and P.~G.~Mehta are with the Coordinated
  Science Laboratory and the Department of Mechanical Science and
  Engineering at the University of Illinois at Urbana-Champaign
  (UIUC); S.~P.~Meyn is with the Department of Electrical and Computer
  Engineering at the University of Florida at Gainesville; Corresponding email: mehtapg@illinois.edu.}
}

\begin{document}

\maketitle
\thispagestyle{empty}
\pagestyle{empty}

\begin{abstract}

Duality between estimation and optimal control is a problem of rich
historical significance.  The first duality principle appears in the
seminal paper of Kalman-Bucy, where the problem of minimum variance
estimation is shown to be dual to a linear quadratic (LQ) optimal
control problem.  Duality offers a constructive proof technique to
derive the Kalman filter equation from the optimal control solution.  
  
\smallskip

This paper generalizes the classical duality result of Kalman-Bucy to
the nonlinear filter: The state evolves as a continuous-time Markov
process and the observation is a nonlinear function of state corrupted
by an additive Gaussian noise.  A dual process is introduced as a
backward stochastic differential equation (BSDE).  The process is used to
transform the problem of minimum variance estimation into an optimal
control problem.  Its solution is obtained from an application of the
maximum principle, and subsequently used to derive the
equation of the nonlinear filter.  The classical duality result of
Kalman-Bucy is shown to be a special case.

\end{abstract}

\section{Introduction}

In Kalman's celebrated paper with Bucy, it is shown that the problem
of minimum variance estimation is dual to a deterministic optimal
control problem~\cite{kalman1961}.  Duality offers a constructive proof technique to
derive the Kalman filter equation from the optimal control solution~\cite[Ch.~7]{astrom1970}.    Apart from the formulation's aesthetic appeal to
control theorists and aficionados of variational techniques, the proof
helps explain why, with the time arrow reversed, the covariance
update equation of the Kalman filter is the same as the dynamic Riccati
equation (DRE) of   optimal control.  Given this, two natural
questions are:  (i) What is the dual optimal control problem for the
  nonlinear filter? and (ii) Can the equation of the nonlinear filter
  be derived from the solution of an optimal control problem?  These
  questions are answered in this paper.  

In classical linear Gaussian settings, dual constructions are of the
  following two types \cite[Sec. 7.3]{bensoussan2018estimation}: 
  	(i) minimum variance estimation,  which was first outlined in
  Kalman-Bucy's original paper, 
  	and 
  	(ii) minimum energy estimation
  whose formulation first appears in~\cite{mortensen1968}.  

Given the historical significance of this area, several 
extensions have been considered over the
decades~\cite{mortensen1968,simon1970duality,FlemingMitter82,fleming1997deterministic,goodwin2005}.
Much work has been done on extending and interpreting the duality for minimum
energy estimation as (i) a MAP estimator~\cite{goodwin2005,todorov2008general}; (ii)
through an application of the \textit{log transformation} to transform the Bellman
equation of optimal control into the Zakai equation of
filtering~\cite{FlemingMitter82,fleming1997deterministic,kappen2016adaptive};
or (iii) based upon the variational Kallianpur-Striebel formula~\cite{mitter2003},~\cite[Lemma
2.2.1]{van2006filtering}.  Based on these extensions, the negative
log-posterior has been shown to have an interpretation as an optimal value function.
Such a formulation is used to derive the equations of nonlinear
smoothing in a companion paper on arxiv~\cite{kim2019smoothing}.

It must be said that none of these earlier results are extensions of duality for the minimum variance estimation problem. 
\spm{killed this since it seemed distracting:  as the
duality principle was first 
described in Kalman-Bucy's paper}
It has been noted in prior work that
	 (i) the dual relationship between the DRE of the
LQ optimal control and the covariance update equation of the Kalman
filter is {\em not} consistent with the interpretation of the negative
log-posterior as a value function,
	 and 
	 (ii) some of the linear
algebraic operations, e.g., the use of matrix transpose to define the
dual system, are not applicable to nonlinear
systems~\cite{todorov2006optimal,todorov2008general}.  
For these reasons, the original duality of Kalman-Bucy is widely understood as an LQG artifact
that does not generalize~\cite{todorov2006optimal}.

\textit{The present paper has a single contribution:} generalization
of the original Kalman-Bucy duality theory to nonlinear filtering.
It is an exact extension, in the sense that the dual optimal control
problem has the same minimum variance structure for linear and
nonlinear filtering problems.  Kalman-Bucy's linear Gaussian result is
shown to be a special case.  
Explicit expressions for the
control Lagrangian and the Hamiltonian are described.  These expressions
are expected to be useful to construct approximate algorithms for
filtering via learning techniques that have become popular of late.     

A related but distinct formulation for duality was proposed by the authors
in a recent paper~\cite{kim2018duality}.  The formulation described in this paper is
original and fixes many of the issues (information structure,
stochastic terms in the objective function) in the earlier paper.
Both the objective function and the BSDE constraint described in this
paper are original.  The main analysis technique -- the use of maximum
principle to derive the nonlinear filter -- is also original.


The outline of the remainder of this paper is as follows: The dual optimal control
problem is proposed in \Sec{sec:problem_formulation}.  The dual
problem is described for two cases:  in the first the state space is finite,  and in the second the state is defined as an 
It\^{o}-diffusion.  The solutions for these two cases appears
in~\Sec{sec:main} and \Sec{sec:main2}, respectively.  A martingale characterization of
the solution appears in~\Sec{sec:mg}.  The proofs are contained in the
Appendix.  


\section{Problem Formulation}
\label{sec:problem_formulation}

\newP{Notation} 
For state-space denoted by $\mathbb{S}$, we let $\clB(\mathbb{S})$ denote the Borel $\sigma$-algebra on $\mathbb{S}$,  and ${\cal  P}(\mathbb{S})$ denote the set of probability measures on $\clB(\mathbb{S})$. Any probability measure $\mu\in{\cal P}(\mathbb{S})$ acts
on a Borel-measurable function $y$ according to $\mu(y) := \int_{\mathbb{S}} y(x)
\mu(\ud x)$.

For a filtration $\clZ := \{\clZ_t : 0\le t\le T\}$
and a measurable space $S$, $\lsq{\clZ}{S}$ denotes the space of
$\clZ$-adapted square-integrable processes taking values in $S$.
Likewise, $L^2_{\clZ_T}(S)$ is the $\clZ_T$-measurable square-integrable
random variables taking values in $S$.  
$C^k(\Re^d;S)$ is the space of $k$-times differentiable functions from
$\Re^d$ to $S$,
 and $L^2(\Re^d;S)$ is the space of square-integrable
(with respect to the Lebesgue measure) function from $\Re^d$ to $S$. 

For a matrix, $\tr(\cdot)$ denotes the trace and
$\dv(\cdot)$ denotes the vector of its diagonal entries.
For a vector, $\dv^\dagger(\cdot)$ denotes a diagonal matrix
with diagonal entries given by the vector.  For a function $y\in
C^2(\Re^d)$, $\frac{\partial y}{\partial x}$ is the gradient vector,
\spm{I'm so surprised you call this a gradient.   I would write gradient = $\nabla y = \frac{\partial y}{\partial x}^\transpose$.  Gradients are column vectors, and derivatives are linear approximations of functions (a row vector for a real valued function). }
 and
$\frac{\partial^2 y}{\partial x^2}$ is the Hessian matrix. For a vector-valued function $f \in C^1(\Re^d\,;\Re^d)$, $\operatorname{div}(f)$ is the divergence of $f$.

\subsection{Filtering problem}

Consider a pair of continuous-time stochastic processes $(X,Z)$. 
The state $X=\{X_t:\, 0\le t\le T \}$ is a Markov process that evolves in
the state-space $\mathbb{S}$.  
 The vector-valued observation process
$Z=\{ Z_t :\, 0\le t\le T \}$ is defined according to the following model:
\begin{equation}\label{eq:obs_model}
Z_t = \int_0^t h(X_s) \ud s + W_t
\end{equation}
where $h:\mathbb{S}\rightarrow \mathbb{R}^m$ is the observation
function and $W=\{W_t:\,t \ge 0\}$ is an $m$-dimensional Wiener process (w.p.) with
covariance matrix $R\succ 0$.  The initial distribution for  $X_0$ is
denoted $\pi_0\in {\cal P}(\mathbb{S})$.

The filtering problem is to compute the conditional 
distribution (posterior) of the state $X_t$ given the filtration (time
history of observations) $\clZ_t := \sigma(Z_s,  0\le s \le t)$.  The posterior
distribution at time $t$ is denoted as $\pi_t$.  It is an element of
${\cal P}(\mathbb{S})$.  For any integrable 
function $f:\mathbb{S}\rightarrow \mathbb{R}$, 
\begin{equation}\label{eq:pit_defn}
\pi_t(f) := {\sf E} (f(X_t) \mid \clZ_t)
\end{equation}

It is well-known that $\pi_t(f)$ is the minimum variance estimator of
$f(X_t)$~\cite[Lemma 5.1]{xiong2008introduction}.
The central question of this paper is to formulate the minimum
variance estimation as a dual optimal control problem. 
\spm{changed `problem' to `formulation' -- don't love it, but don't love problems, either!}
 These formulations are described in 
\Sec{sec:euclidean_ss} and \Sec{subsec:problem-MC}, for the Euclidean and the finite
state-space settings, respectively.  The duality
principle is then presented in \Sec{sec:duality}.  The relationship to
the well known linear Gaussian case is discussed in \Sec{sec:LG}.

\subsection{It\^{o} diffusion on the Euclidean space} 
\label{sec:euclidean_ss}

\spm{It is unfortuante that we are swapping the order of the two cases back and forth.  Would be better to do finite, then diffusion}

The state
process $X$ evolves 
on    $\mathbb{S}=\Re^d$
according to the It\^o stochastic differential
equation (SDE) 
\[
\ud X_t = a (X_t) \ud t + \sigma(X_t) \ud B_t, \quad X_0\sim \pi_0
\]
where $B=\{B_t:\,t \ge 0\}$ is a vector valued standard w.p., 
and $a(\cdot),\;\sigma(\cdot)$ are $C^2$ functions of appropriate
dimensions.  It is assumed
that $W$, $B$,  $X_0$ are mutually independent.  

\spm{Moved this to notation: 
Any probability measure $\pi\in{\cal  P}(\mathbb{S})$ acts
on a Borel-measurable function $y$ according to $\pi(y) := \int_{\Re^d} y(x)
\pi(\ud x)$.
}


The differential generator of $X$, denoted as $\clA$,
acts on $C^2$ functions in its domain according to 
$$
(\clA y)(x):= a^\top(x) \frac{\partial y}{\partial
	x} (x) + \half \tr\Big(\sigma(x)\sigma^\top(x)\frac{\partial^2y}{\partial x^2} (x)\Big)
$$
It is
assumed that $\clA$ is an elliptic operator:  there is $\varepsilon>0 $ such that $\sigma(x)\sigma^\top(x)\ge \varepsilon I$ for all $x\in\Re^d$.
\spm{This is what you want, right?}

For functions $y\in C^1(\Re^d;\Re)$, $v\in C (\Re^d;\Re^{m})$, 
$u\in\Re^m$, and $x\in\Re^d$, a \textit{cost function} is defined as follows:
\begin{equation}\label{eq:lagrangian-euclide}
\ell(y,v,u\,;x) = \half \Big|\sigma^\top(x) \frac{\partial y}{\partial
    x}(x) \Big|^2 + \half (u+v(x))^\top R (u+ v(x))
\end{equation}

\newP{Dual optimal control problem}  
\begin{subequations}\label{eq:opt-cont-euclide}
	\begin{align}
	&\mathop{\text{Min}}_{U\in\;\clU}\ \sJ(U) =
          \E\Big(\half |Y_0(X_0)-\pi_0(Y_0)|^2 + \int_0^T \ell (Y_t,V_t,U_t\,;X_t) \ud t \Big)\label{eq:opt-cont-euclide-a}\\
	&\text{Subj.}\ \ud Y_t(x)  = - \big((\clA Y_t)(x) + h^\top(x) (U_t + V_t(x))\big)\ud t +V_t^\top(x)\ud Z_t \nonumber \\
	&\quad\quad\;\;\;\; Y_T (x) = f(x),\quad \forall\, x\in\Re^d
\label{eq:opt-cont-euclide-b}
	\end{align}
\end{subequations}
The constraint
\eqref{eq:opt-cont-euclide-b}
 is a
backward stochastic partial differential equation (BSPDE) with
boundary condition prescribed at the
terminal time $T$.    The function $f$ appearing in this boundary condition is allowed to be random, with $f\in L^2_{\clZ_T}(L^2(\Re^d;\Re))$.
\spm{Attempt at emphasis of random $f$ here}

   The admissible set of control input is as follows:
\[
{\cal U} := \lsq{\clZ}{\Re^m}
\]
The solution $(Y,V):=\{(Y_t(x),V_t(x))\,:\, 0\le t \le T,\;x\in\Re^d\}$ of the BPSDE is
adapted to the filtration $\clZ$.  It is 
an element of $L^2_{\clZ}([0,T];L^2(\Re^d;\Re))\times
L^2_{\clZ}([0,T];L^2(\Re^d;\Re^m))$; cf.,~\cite{ma1997adapted}.   




\subsection{Finite state-space}
\label{subsec:problem-MC}

The continuous-time process evolves on the finite state space $\mathbb{S} = \{e_1,\ldots,e_d\}$.  Its statistics are characterized by the  initial distribution $\pi_0\in {\cal P}(\mathbb{S})$  and the
row-stochastic rate matrix 
$A\in\Re^{d\times d}$. 


The dual of ${\cal  P}(\mathbb{S})$ is the space of all  functions on $\mathbb{S}$,  which can be identified with $\Re^d$:    Any function $y:   \mathbb{S}\to \Re$ is determined by its values at the basis vectors $\{e_i\}$, and for any $\pi \in {\cal  P}(\mathbb{S})$, the expectation can be expressed as a dot product:  $\pi(y) = \sum \pi(e_i) y(e_i)$.    Similarly, the observation function can be expressed $h(x) = H^\top x$,  $x\in \bS$,  where
$H\in\Re^{d\times m}$.   
\spm{I understood the intention of the previous essay, but found it hard to read.  Feel free to restore text if you feel the revision is worse!}
\jin{rearranging and minor edit}



We also define a $d\times d$ matrix for each $i$ as follows:
\begin{equation}\label{eq:Q-explicit}
Q(e_i) := \sum_{j=1}^dA_{ij}(e_j-e_i)(e_j-e_i)^\top,\quad i = 1,\ldots,d
\end{equation}
The cost function in this case is defined as
\begin{equation}
\ell(y,v,u\,;x) = \half y^\top Q(x) y + \half (u+v^\top x)^\top R (u+v^\top x)
\end{equation}
where $y\in\Re^d$, $v \in \Re^{d\times m}$, $u\in\Re^m$, and $x\in \mathbb{S}\subset
  \Re^d$.

\newP{Dual optimal control problem}  
\begin{subequations}\label{eq:opt-cont-finite}
	\begin{align}
	&\mathop{\text{Min}}_{U\in\,\clU}\ \sJ(U) =
          \E\Big(\half|Y_0^\top(X_0-\pi_0)|^2 + \int_0^T \ell (Y_t,V_t,U_t\,;X_t) \ud t \Big)\label{eq:opt-cont-finite-a}\\
	&\text{Subj. }\  \ud Y_t = -\big(AY_t+HU_t + \dv(HV_t^\top)\big) \ud t + V_t \ud Z_t,\quad Y_T = f \label{eq:opt-cont-finite-b}
	\end{align}
\end{subequations}
The constraint is a backward stochastic
differential equation (BSDE) with terminal condition as before. 
We write $Y_T=f\in L^2_{\clZ_T}(\Re^d)$,  with our convention that functions on $\bS$ are identified with vectors in $\Re^d$.

The admissible set of control input is $\clU=\lsq{\clZ}{\Re^m}$ and
the solution pair 
$\{(Y_t,V_t)\,:\, 0\le t \le T\}=:(Y,V)\in
\lsq{\clZ}{\Re^d}\times\lsq{\clZ}{\Re^{d\times m}}$;
cf.~\cite[Ch. 7]{yong1999stochastic}.  


\medskip

\subsection{Duality relationship} 
\label{sec:duality}

Suppose $Z$ is defined according the observation
model~\eqref{eq:obs_model}.  Consider an admissible control input
$U\in {\cal U}$ and define $Y_0$ via solution of the BSDE -- Eq.~\eqref{eq:opt-cont-euclide-b} for the Euclidean case,
 and
Eq.~\eqref{eq:opt-cont-finite-b} for the finite case.  Assume the
following linear structure of the estimator:
\begin{equation}\label{eq:opt_est}
S_T = \pi_0(Y_0) - \int_0^T U_t^\top \ud Z_t
\end{equation}    

The duality relationship is expressed in the following proposition
whose proof appears in \Appendix{apdx:opt_control}:

\begin{proposition}\label{prop:minimum_variance}
Consider the observation model~\eqref{eq:obs_model}, the linear
estimator~\eqref{eq:opt_est}, together with the dual optimal control
problem. 
Then for any
choice of admissible control $U\in {\cal U}$: 
\begin{equation*}
\sJ(U) = \half \E (|S_T - f(X_T)|^2)
\end{equation*}	
\end{proposition}

\medskip

Thus, formally, the problem of obtaining the minimum variance estimate $S_T$
of $f(X_T)$ (minimizer of the right-hand side of the equality) is 
converted into the problem of finding the optimal control $U$
(minimizer of the left-hand side of the identity).  However, there is a
subtle problem:  It is not apriori clear whether there exists a
$U\in{\cal U}$ such that $S_T = \pi_T(f)$\footnote{This will be true, e.g.,
if {\em all} $\clZ_T$-measurable random variable have a representation of the
form~\eqref{eq:opt_est}.}. 
\spm{I don't like burying this as a footnote, but I'm out of time to edit this carefully}
 In this paper, the following
assumption is made: 


\newP{Assumption A1}
For each fixed terminal time $T>0$ and function $f\in L^2_{{\clZ}_T}$, there exists a $U \in {\cal U}$ such
that $S_T =  \pi_T(f)$.

\medskip

Under this assumption, the following is proved in 
\Appendix{apdx:proof-prop2}: 





\begin{proposition}\label{prop:optimal-solution}
        Consider the dual optimal control problem. Suppose
$U=\{U_t:0\le t \le T\}$ is the optimal control input and that $Y =
\{Y_t: 0\le t \le T\}$ is the associated optimal trajectory
obtained as a solution of the BSDE.  Then for all $t \in [0,T]$:
	\begin{equation}\label{eq:estimator-t}
	\pi_t(Y_t) = \pi_0(Y_0) - \int_0^t U_s^\top\ud Z_s 
	\end{equation}
\end{proposition}

\medskip



	

\subsection{Linear-Gaussian case} 
\label{sec:LG}

The linear-Gaussian case assumes the following model:
\begin{enumerate}
	\item The drift in the It\^{o} diffusion is linear in $x$.  That is, 
\[ 
a(x)=A x \;\; \text{and}\;\; h(x) = H x
\] 
where $A\in\Re^{d\times d}$ and $H\in\Re^{m\times d}$.
	\item The coefficient of the process noise is a constant matrix,  $
\sigma(x) \equiv \sigma$.  We denote $Q:=\sigma\sigma^\top\in\Re^{d\times
          d}$.
	\item The prior $\pi_0$ is a Gaussian distribution with mean
          $m_0\in\Re^{d}$ and variance $\Sigma_0\in\Re^{d\times d}$.
\end{enumerate}

Classical Kalman-Bucy duality is concerned with the problem of
constructing a minimum variance estimator for the random variable
$f^\top X_T$ where $f\in \Re^d$ is a given deterministic
vector~\cite{kalman1961}.  Therefore, we also have
\begin{enumerate}
	\setcounter{enumi}{3}
	\item The terminal condition in~\eqref{eq:opt-cont-euclide-b} is a linear
	function $f^\top x$.  
\end{enumerate}

\medskip

We impose the following restrictions:
\begin{enumerate}
\item The control input $U=u$ is restricted to be a deterministic
  function of time (in particular, it does not depend upon the
  observations).  Such a control is trivially $\clZ$-adapted hence
  admissible.  
  For such a control
  input, with deterministic $f$, the solution $Y=y$ of the BSPDE is a
  deterministic function of time, and $V=0$.  The BSPDE becomes a
  PDE:
\begin{equation}\label{eq:det_pde}
\frac{\partial y_t}{\partial t}(x) = -(\clA y_t)(x) -  h^\top(x)
u_t,\;\; y_T(x)\equiv f^\top x \;\; \forall x\in\Re^d
\end{equation}
where the lower-case notation is used to stress the fact that $u$
and $y$ are
now deterministic functions of time.  
\item Without loss of generality, it suffices to consider the
  restriction of the optimal control
  problem~\eqref{eq:opt-cont-euclide} on a finite ($d-$) dimensional
  subspace of the function space:
\[
{\cal V}:=\{\tilde{y} \;:\;  \tilde{y}(x)=y^\top x \;\; \forall\, x\in\Re^d\;\text{where}\;y\in\Re^d\}
\]  
It is easy to see that ${\cal V}$ is an invariant subspace for the
dynamics~\eqref{eq:det_pde}.  On ${\cal V}$, the PDE reduces to an ODE:  
\[
\frac{\ud y_t}{\ud t}
= -A^\top y_t - H^\top u_t,\quad y_T = f 
\]
and the cost function becomes
\[
{\sf L} (y,u) := \half y^\top Q y + \half u^\top R u
\]
It no longer depends upon $x$.
\end{enumerate}




In summary, the  optimal control problem~\eqref{eq:opt-cont-euclide}
reduces to the deterministic LQ problem of classical duality:
\begin{align*}
&\mathop{\text{Minimize}}_{u}:\quad \sJ(u) = \half \ y_0^\top \Sigma_0 y_0 + \int_0^{T} {\sf L}(y_t,u_t)\ud t \\
&\text{Subject to}:\quad  \frac{\ud y_t}{\ud t}
= -A^\top y_t - H^\top u_t,\quad y_T = f 
\end{align*}

The solution of
the optimal control problem yields the optimal control input $u$, 
along with the vector $y_0$ that determines the minimum-variance
estimator:
\begin{align*}
S_T &= \pi_0(y_0^\top x) - \int_0^T u_t^{\top} \ud Z_t
= y_0^\top m_0 - \int_0^T u_t^{\top} \ud Z_t
\end{align*}
The Kalman filter is obtained by expressing $\{S_t(f) : t\ge 0,\ f\in\Re^d\}$ as the solution to a linear SDE~\cite[Ch.~7]{astrom1970}.      

\medskip

\begin{remark}
Consider the dual optimal control problem~\eqref{eq:opt-cont-finite}
for the finite state space model.  Suppose one only allows
control inputs $U$ that are deterministic functions of time.  
In this case, $Y$ is a deterministic function
of time and $V=0$.  Consequently, the objective function
in~\eqref{eq:opt-cont-finite-a} simplifies
\spm{sometimes you really want square brackets -- here, all of the ))) looked odd}
\begin{equation}\label{eq:DLQ-finite}
\sJ(U) = \half Y_0^\top \Sigma_0Y_0 + \int_0^T \half U_t^\top R U_t +
\half Y_t^\top \E ( Q(X_t) ) Y_t \ud t
\end{equation}
where $\Sigma_0 := \E((X_0-\pi_0)(X_0-\pi_0)^\top)$ and $\E (Q(X))$ is
the quadratic variation process for $X$. 
The resulting problem is a deterministic LQ problem whose optimal solution
$\{U_t:0\le t \le T\}$ will (in general) yield a sub-optimal estimate $S_T$
using~\eqref{eq:opt_est}.  
In Appendix~\ref{apdx:KF}, the solution is used to derive a Kalman filter for the Markov chain.  Such sub-optimal
filters for Markov chains have been applied in~\cite{lipkrirub84,chebusmey17a}.
\end{remark}
\medskip

We  have now set the stage to derive the nonlinear filter via the
solution to the dual optimal control problem.  We describe the
solution for the finite state case first in \Sec{sec:main}.  Although technically and
notationally more challenging, the considerations for the
Euclidean case are entirely analogous and described in~\Sec{sec:main2}.

\section{Solution for the finite-state case}\label{sec:main}

\subsection{Standard form of the optimal control problem}

Consider the dual optimal control problem~\eqref{eq:opt-cont-finite} associated with the finite state space model.  
There is only one natural filtration in this problem,  the
filtration $\clZ$  generated by the observation process.  
\spm{Most readers would assume that this is the only natural filtration -- I don't understand the comment}
The `state' of the
optimal control problem $(Y,V)$ is adapted to the filtration by
construction.  So, the problem is fully observed.  However, the
problem is not in its standard form \cite[Def. 5.4]{pardoux2014stochastic}.  
\spm{This text confuses me, and I don't think it is essential:  as the form pertains to the optimal
control of BSDEs~}
There are two issues:
\begin{enumerate}
\item The stochastic process $Z$ on the right-hand side of the BSDE~\eqref{eq:opt-cont-finite-b} is not a Wiener process.
\item The cost function in~\eqref{eq:opt-cont-finite-b} depends upon the exogenous process $X$ which is not adapted to $\clZ$.   
\end{enumerate}  

In order to resolve these issues and express the optimal control
problem in a standard form, we introduce three stochastic processes
$\pi:=\{\pi_t\in{\cal P}(\mathbb{S})\subset\Re^d:0\le t\le T\}$,
$\Sigma:=\{\Sigma_t\in\Re^{d\times d} : 0\le t\le T\}$, and
$\innov:=\{\innov_t\in\Re^m:0\le t\le T\}$ as follows:
\begin{subequations}\label{eq:pi_and_I}
	\begin{flalign}
	&\text{filter}:  &\pi_t &:= E(X_t|\clZ_t) &\label{eq:
          pi_and_I-a}\\
        &\text{covariance}: &\Sigma_t&:= \dv^\dagger(\pi_t)-\pi_t\pi_t^\top& \label{eq:pi_and_I-b}\\
	&\text{innovation}: &\innov_t &:= Z_t - \int_0^t \pi_s(h) \ud s & \label{eq:pi_and_I-c}
	\end{flalign}
\end{subequations}
From the standard filtering theory, it is well known that i) $\innov$ is a
Wiener process that is adapted to $\clZ$~\cite[Lemma~5.6]{xiong2008introduction},
\spm{removed "with respect" and replaced by adapted.  OK? }
\jin{Okay!}
and ii) the filtration
$\sigma(\innov_s: 0\le s\le t)$ generated by the innovation process 
equals $\clZ_t$ for all $t\in[0,T]$~\cite{mitter1982nonlinear}.   

We therefore express the BSDE constraint as
\begin{align*}
\ud Y_t &= -\big(AY_t+HU_t + \dv(HV_t^\top)\big) \ud t + V_t \ud Z_t \\
&=-\big(AY_t+HU_t + \dv(HV_t^\top)-V_tH^\top \pi_t\big) \ud t + V_t \ud \innov_t
\end{align*}
The solution $(Y,V)$ is adapted to $\clZ$, now interpreted as the filtration generated by the innovation process.  

In order to remove the explicit dependence of the cost function on the
non-adapted process $X$, we use the tower property of the conditional
expectation:
\[
\sJ(U) = \half \E (|Y_0(X_0) -\pi_0(Y_0)|^2) + \int_0^T \E( \ell (Y_t,V_t,U_t\,;X_t)|\clZ_t) \ud t 
\]
Denote ${\cal{L}}
(Y_t,V_t,U_t\,;\pi_t) \eqdef \E( \ell (Y_t,V_t,U_t\,;X_t)|\clZ_t)$.

Because $Y,V,U$ are all $\clZ$-adapted, it is a straightforward calculation to
see that  
\[
{\cal{L}} (y,v,u\,;\mu) = \half y^\top \mu(Q) y + \half u^\top R u+ u^\top Rv^\top \mu + \half \mu^\top \dv(v R v^\top)
\]
where $\mu(Q) := \dv^\dagger (A^\top\mu) - A^\top
\dv^\dagger(\mu)-\dv^\dagger(\mu) A$.  The function ${\cal{L}}$
is the control Lagrangian.  

We are now ready to state the standard form of the optimal control
problem for the finite state-space  case:

\newP{Dual optimal control problem (standard form)}  
\begin{subequations}\label{eq:opt-cont-pi}
\begin{align}
&\mathop{\text{Min}}_{U\in\clU}\ \sJ(U) = \E\Big(\half
  Y_0^\top\Sigma_0Y_0 + \int_0^T {\cal L}(Y_t,V_t,U_t;\pi_t) \ud t \Big)\label{eq:opt-cont-pi-a}\\
&\text{Subj.}\  \ud Y_t = -\big(AY_t+HU_t + \dv(HV_t^\top)-V_tH^\top \pi_t\big) \ud t + V_t \ud \innov_t, \nonumber\\
&\quad \quad \;\;\;\; Y_T = f \label{eq:opt-cont-pi-b}
\end{align}
\end{subequations}

In its standard form, all the processes are adapted to the filtration
$\clZ$ and moreover the stochastic process $\innov$ on the right-hand side
of the BSDE is a Wiener process
with respect to this filtration.

\subsection{Solution using the maximum principle} 

The Hamiltonian ${\cal H}:\Re^d\times \Re^{d\times m}\times \Re^m \times \Re^d \times {\cal P}(\mathbb{S})\to \Re$ is defined as follows:
\begin{align*}
{\cal H}(y,v&,u,p\,; \mu) \\
&= p^\top\big(-Ay - Hu - \dv(Hv^\top) + vH^\top \mu\big) - {\cal L}(y,v,u\,; \mu)
\end{align*}

A characterization of the optimal input is contained in the following.  Its proof, based on an  application of the maximum
principle~\cite{peng1993backward}, appears in the \Appendix{apdx:maximum}.  

\medskip

\begin{theorem}
\label{thm:optimal-solution} 

Consider the optimal control problem~\eqref{eq:opt-cont-pi}.  Suppose
$U=\{U_t:0\le t \le T\}$ is the optimal control input and that $(Y,V) =
\{(Y_t,V_t): 0\le t \le T\}$ is the associated optimal trajectory
obtained as a solution of the BSDE~\eqref{eq:opt-cont-pi-b}.  Then
there exists a $\clZ$-adapted vector valued process $P=\{P_t:0\le t \le T\}$ such
that   
\begin{equation}\label{eq:opt-cont-soln}
U_t = -R^{-1}H^\top P_t - V_t^\top \pi_t
\end{equation}
where      $P$  and  $Y$  satisfy
\spm{Note minor rearrangements.
\\
Why not give $P$ a name (co-state)? And what about poor $Y$?   }
\begin{subequations}
\label{eq:Hamilton_eqns}
\begin{flalign}
&\text{(forward)}  & \ud P_t &= -{\cal H}_y(Y_t, V_t,U_t,P_t\,;\pi_t)\ud t \nonumber \\
& & &\ \  - {\cal H}_v(Y_t, V_t,U_t,P_t\,;\pi_t) R^{-1}\ud \innov_t
& \label{eq:Hamilton_eqns-a}\\
&\text{(backward)}  &\ud Y_t &= {\cal H}_p(Y_t,V_t,U_t,P_t\,;\pi_t) \ud t + V_t\ud \innov_t & \label{eq:Hamilton_eqns-b}\\
&\text{(boundary)} & P_0 &= \Sigma_0 Y_0 , \quad Y_T = f \label{eq:Hamilton_eqns-c}&
\end{flalign}
\end{subequations}
\end{theorem}

\medskip

\begin{remark} \label{rm:linear-structure}
From linear optimal control theory, it is known that $P_t = M_t Y_t$
(see~\cite[Sec. 6.6]{yong1999stochastic}) where $M:=\{M_t\in\Re^{d\times d}:0\le t\le T\}$ is a $\clZ$-adapted
matrix-valued process.  The boundary condition $P_0=\Sigma_0 Y_0$
suggests that $M=\Sigma$.  This is indeed the case as shown in the
proof of Theorem~\ref{thm:Wonham},
where the following equation is  
derived (see~\eqref{eq:DRE-for-NLF}):
\begin{align*}
\ud \Sigma_t
=&\big(A^\top\Sigma_t + \Sigma_t A + \pi_t(Q) -\Sigma_tHR^{-1}H^\top\Sigma_t\big)\ud t\\
&+\dv^\dagger(\Sigma_t H R^{-1} \ud \innov_t) - \Sigma_t H R^{-1} \ud \innov_t \pi_t^\top- \pi_t\ud \innov_t^\top R^{-1}H^\top\Sigma_t
\end{align*}
This is the DRE of the nonlinear filter.

\end{remark}




\subsection{Derivation of the nonlinear filter} 

From \Prop{prop:optimal-solution}, using the formula~\eqref{eq:opt-cont-soln} for the optimal control,
\begin{equation}\label{eq:pi_t_y_t_est}
\pi_t^\top Y_t = \pi_0^\top Y_0 + \int_0^t ( P_s^\top H R^{-1} + \pi_s^\top V_s) \ud Z_s
\end{equation}
This formula is used to derive the Wonham filter~\cite{wonham1964}.  The proof of the
following theorem appears in the \Appendix{apdx:Wonham}.


\spm{It is "Hamilton's equation", without 'the',  or "the Hamilton equation"}

\begin{theorem}[Nonlinear filter]
\label{thm:Wonham} 
Consider the optimal estimator~\eqref{eq:pi_t_y_t_est} where $(Y,V)$
and $P$ solve   Hamilton's equations~\eqref{eq:Hamilton_eqns}.  Then  
\begin{equation}
\label{eq:Wonham}
\ud\pi_t = A^\top\pi_t\ud t + \Sigma_tHR^{-1}\ud \innov_t
\end{equation}
and furthermore $P_t=\Sigma_t Y_t$ for all $t\in[0,T]$.  
\end{theorem}






\section{Solution for the Euclidean case}\label{sec:main2}

As in the finite state-space case, the starting point is the
definition of the Lagrangian: ${\cal{L}}
(Y_t,V_t,U_t\,;\pi_t):=\E( \ell (Y_t,V_t,U_t\,;X_t)|\clZ_t)$ where the
cost function $\ell$ for the Euclidean case is defined
in~\eqref{eq:lagrangian-euclide} and $\pi_t$ is the conditional
distribution at time $t$ defined according to~\eqref{eq:pit_defn}.  Explicitly,  
\begin{align*}
{\cal L} & (y,v,u\,;\mu) \\
& =\half \int_{\Re^d} \left( \left|\sigma^\top(x) \frac{\partial y}{\partial
  x}(x)\right|^2 + (u+v(x))^\top R (u+ v(x)) \right) \mu (\ud x)
\end{align*}

The standard form of the optimal control problem is as follows:

\newP{Dual optimal control problem (standard form)}
\begin{subequations}\label{eq:opt-cont-euc-pi}
	\begin{align}
	&\mathop{\text{Min}}_{U\in\clU}\ \sJ(U) = \E\Big(\half
	\pi_0(|Y_0-\pi_0(Y_0)|^2) + \int_0^T {\cal L}(Y_t,V_t,U_t;\pi_t) \ud t \Big)\label{eq:opt-cont-euc-pi-a}\\
	&\text{Subj.}\  \ud Y_t = -\big((\clA Y_t)(x) + h^\top(x)(U_t + V_t(x)) + V_t^\top(x)\pi_t(h)\big) \ud t \nonumber\\
	&\quad\quad\quad\quad\quad + V_t^\top(x) \ud \innov_t,\quad Y_T(x) = f(x)\quad \forall\, x\in\Re^d \label{eq:opt-cont-euc-pi-b}
	\end{align}
\end{subequations}

\newP{Assumption A2} The (generic) measure $\mu$ is absolutely
continuous with respect to   Lebesgue measure.

\newP{Notation} The Radon-Nikodyn derivative is denoted as
$\tilde{\mu}(x) := \frac{\ud \mu}{\ud x}(x)
$. 
Consequently, $\mu(f) = \int_{\Re^d} f(x)
\tilde{\mu}(x) \ud x =: \langle \tilde{\mu},f \rangle$.  In the remainder of this section, with a
slight abuse of notation, we will drop the tilde to simply write $\langle
\mu,f \rangle$. 

\medskip

The co-state $p\in L^2(\Re^d;\Re)$ and the Hamiltonian are defined as follows:
\spm{now you give $p$ a name -- why not at the start?}
\begin{align*}
{\cal H}(y,v&,u,p\,;\mu) = \langle p,-\clA y - h^\top (u + v) + v^\top \langle\mu,h\rangle)\rangle- {\cal L}(y,v,u\,;\mu)
\end{align*}
Hamilton's equations are described in the following, while the proof appears in \Appendix{apdx:L2-derivatives}:





\begin{theorem}\label{thm:optimal-solution-euc}
Consider the optimal control problem~\eqref{eq:opt-cont-euc-pi}.  Suppose
$U=\{U_t:0\le t \le T\}$ is the optimal control input and the $(Y,V) =
\{(Y_t,V_t): 0\le t \le T\}$ is the associated optimal solution
obtained by solving BSPDE~\eqref{eq:opt-cont-euc-pi-b}.
Then there exists a $\clZ$-adapted function-valued process $P=\{P_t:0\le t \le T\}$ such
that
\begin{subequations}\label{eq:Hamilton_eqns_euc}
	\begin{flalign}
	&\text{(forward)}  & \ud P_t &= -{\cal H}_y(Y_t, V_t,U_t,P_t\,;\pi_t)\ud t \nonumber \\
	& & &\ - {\cal H}_v^\top(Y_t, V_t,U_t,P_t\,;\pi_t) R^{-1}\ud \innov_t
	& \label{eq:Hamilton_eqns_euc-a}\\
	&\text{(backward)}  &\ud Y_t &= {\cal H}_p(Y_t,V_t,U_t,P_t\,;\pi_t) \ud t + V_t\ud \innov_t & \label{eq:Hamilton_eqns_euc-b}\\
	&\text{(boundary)} & P_0(x&) = \pi_0(x)\big(Y_0(x)-\langle\pi_0,Y_0\rangle\big)& \nonumber\\
	&  & Y_T(x&) = f(x)\quad \forall \, x\in\Re^d & \label{eq:Hamilton_eqns_euc-c}
	\end{flalign}
\end{subequations}
where the optimal control is given by
\begin{equation}\label{eq:opt-cont-soln-euc}
U_t = -R^{-1}\langle P_t,h\rangle - \langle \pi_t, V_t \rangle 
\end{equation}
\end{theorem}

\medskip
\medskip


Using the result of Prop.~\ref{prop:optimal-solution}, the optimal
estimator is
\begin{equation}\label{eq:pi_t-y_t-euc}
\langle\pi_t,Y_t\rangle = \langle \pi_0,Y_0\rangle + \int_0^t \big(R^{-1}\langle P_s, h \rangle + \langle \pi_s, V_s \rangle \big)^\top \ud Z_s
\end{equation}
As before, the formula is used to derive the equation for the Kushner filter equation~\cite{kushner1967dynamical}. The proof appears in \Appendix{apdx:Kushner}:
\begin{theorem}\label{thm:Kushner}
Consider the optimal estimator~\eqref{eq:pi_t-y_t-euc} where $(Y,V)$
and $P$ solve Hamilton's
equation~\eqref{eq:Hamilton_eqns_euc}. Then the conditional density
solves the SPDE:
$$
\ud \pi_t(x) = (A^\dagger\pi_t)(x)\ud t + \pi_t(x)(h(x)-\langle \pi_t,h\rangle)R^{-1}\ud \innov_t
$$
and furthermore $P_t(x) = \pi_t(x) \big(Y_t(x) - \langle \pi_t, Y_t
\rangle\big)$ for all $x\in\Re^d$ and $t\in[0,T]$.
\end{theorem}

\medskip


\section{A martingale characterization}
\label{sec:mg}

For the finite state-space case, define the function ${\cal V}:\Re^d\times
{\cal P}(\Re^d)\rightarrow \Re$ as follows:
\[
{\cal V}(y\,;\mu) = \half \sum_{i=1}^d |y(i) - \mu^\top y|^2 \mu(i)
\]
For the Euclidean case, the analogous function is as follows:
\[
{\cal V}(y\,;\mu) = \half \int_{\Re^d} |y(x)-\langle \mu, y \rangle |^2
\mu(x) \ud x
\]

In the statement of the following theorem, $U^{*}$ denotes the optimal
control input as defined according to the
formula~\eqref{eq:opt-cont-soln} for the finite-state-space  and the
formula~\eqref{eq:opt-cont-soln-euc} for the Euclidean case.  The
proof appears in the \Appendix{apdx:martingale}.

\medskip

\spm{This is a beautiful result.   I wish it were a bit more visible.}
   
\begin{theorem}
\label{thm:martingale-DP}

Suppose $\pi$ is the posterior process and $(Y,V)$ is the $\clZ$-adapted solution of the dual BSDE.  
For every input $U\in{\cal  U}$, the process
\[
{\cal M}_t^U = {\cal V}(Y_t\,;\pi_t) - \int_0^t {\cal L}(Y_s,V_s,U_s\,;\pi_s)\ud s
\]
is a supermartingale; it is a martingale if and only if $U=U^*$.
Consequently, 
\[
\sJ(U) \geq \E\big({\cal V}(f\,;\pi_T)\big)
\]
with equality if and only if $U=U^*$. 
Consequently,  the right-hand side is the value function for the dual optimal control problem.
\end{theorem}

\bibliographystyle{IEEEtran}
\bibliography{duality,backward_sde}


\appendix

\section{Appendix}

\subsection{Proof of Prop.~\ref{prop:minimum_variance}}\label{apdx:opt_control}

\newP{Euclidean case} For a random function $Y_t(\,\cdot\,)$,
the
It\^{o}-Wentzell theorem~\cite[Thm. 1.17]{rozovskiui2018stochastic}
gives the formula for the differential $\ud (Y_t(X_t))$
\begin{align*}
= & -\big(\clA Y_t(X_t) + U_t^\top h(X_t) + V_t^\top h(X_t)\big)\ud t + V_t(X_t) \ud Z_t\\
&\quad+ a^\top(X_t)\frac{\partial Y_t}{\partial x}(X_t)\ud t + \big(\sigma^\top(X_t)\frac{\partial Y_t}{\partial x}(X_t)\big)^\top \ud B_t \\
&\quad + \half \tr\big(\sigma\sigma^\top(X_t)\frac{\partial^2
  Y_t}{\partial x^2}(X_t)\big)\ud t 
+ \tr\big(\frac{\partial V_t}{\partial x}\sigma(X_t) \ud [ B, W ]_t\big)\\
=& -U_t^\top \ud Z_t + (U_t + V_t(X_t))^\top \ud W_t + \big(\sigma^\top\frac{\partial Y_t}{\partial x}(X_t)\big)^\top \ud B_t
\end{align*}
Integrating over $[0,T]$,
\begin{align*}
f(X_T) &= Y_0(X_0) - \int_0^T U_t^\top \ud Z_t 
\\
&\qquad 
	+\int_0^T (U_t + V_t(X_t))^\top \ud W_t + \big(\sigma^\top\frac{\partial Y_t}{\partial x}(X_t)\big)^\top \ud B_t
\end{align*}
Using the formula~\eqref{eq:opt_est} for the estimator, the error 
\begin{align*}
f(X_T) - S_T &= Y_0(X_0) - \pi_0(Y_0) + \int_0^T \big(\sigma^\top\frac{\partial Y_t}{\partial x}(X_t)\big)^\top \ud B_t\\
&\qquad
	+\int_0^T(U_t + V_t(X_t))^\top \ud W_t
\end{align*}
Upon squaring and taking an expectation,
\begin{align*}
 \E\big(|S_T&-f(X_T)|^2\big) \\
 &= \E\Big(|y_0(X_0)-\pi_0(y_0)|^2 + \int_0^T\big|\sigma^\top(X_t)\frac{\partial Y_t}{\partial
  x}(X_t)\big|^2\ud t \\ 
& \qquad \quad +\int_0^T (U_t + V_t^\top(X_t))^\top R (U_t + V_t^\top(X_t)) \ud t\Big)
\end{align*}

\newP{Finite state-space  case} 
A martingale $N = \{N_t: t\ge 0\}$ is defined as follows:
$$
N_t = X_t - \int_0^t A^\top X_s\ud s
$$
whose quadratic variation is
$
[ N ]_t = \int_0^t Q(X_s)\ud s
$  (see  \eqref{eq:Q-explicit} for the definition of $Q$). 
   It\^{o}'s product formula gives
\begin{align*}
\ud (Y_t^\top X_t) 
&=-U_t^\top \ud Z_t + U_t^\top \ud W_t + X_t^\top V_t \ud W_t + Y_t^\top \ud N_t
\end{align*}
Using this, together with~\eqref{eq:opt_est} for the estimator,
results in the following error equation: 
$$
f^\top X_T - S_T = Y_0^\top(X_0-\pi_0) + \int_0^T (U_t + V_t^\top X_t)^\top \ud W_t + Y_t^\top \ud N_t
$$
Upon squaring and taking an expectation,
\begin{align*}
\E\big[ |f^\top X_T &- S_T|^2\big] 
= \E\Big[ \big(Y_0^\top(X_0-\pi_0)\big)^2 \\
& + \int_0^T \Big((U_t + V_t^\top X_t)^\top R (U_t + V_t^\top X_t) + Y_t^\top Q(X_t) Y_t\Big) \ud t \Big]
\end{align*}

\subsection{Proof of Prop.~\ref{prop:optimal-solution}}\label{apdx:proof-prop2}


Suppose $U\in {\cal U}$.  Define
\[
S_t:=\pi_0(Y_0) - \int_0^t U_s^\top\ud Z_s 
\]
Using~\Prop{prop:minimum_variance},
\begin{align*}
\half  \E(|S_t&-Y_t(X_t)|^2) = \\
& \quad \E\Big(\half|Y_0(X_0)-\pi_0(Y_0)|^2 + \int_0^t \ell (Y_s,V_s,U_s\,;X_s)\ud s\Big) 
\end{align*}

By the dynamic programming principle, if $\{U_s: 0 \le s \le T\}$ is
an optimal control input over the time-horizon $[0,T]$ then $\{U_s: 0 \le s \le t\}$
minimizes the right-hand side over all $\clZ$-adapted control inputs.  

It then follows from Assumption~(A1) that $S_t=\pi_t(Y_t)$.  



\subsection{Derivation of the Kalman filter for the Markov process}
\label{apdx:KF}

\def\bSig{\bar{\Sigma}}
\def\barX{\bar{X}}
\def\barY{\bar{Y}}   
\def\barS{\bar{S}}

The optimal control solution is given by
\begin{align*}
&\text{(DRE)} \; \frac{\ud}{\ud t}\bSig_t  = \bSig_t A + A^\top \bSig_t -
                                        \bSig_t HR^{-1} H^\top \bSig_t
                                        + \E(Q(X_t)),\;\; \bSig_0 =
                                        \Sigma_0\\
& \text{(Opt. control)}  \quad U_t = -R^{-1} H^\top \bSig_t Y_t
\end{align*} 
Upon substituting the 
solution into the estimator~\eqref{eq:opt_est}
\[
S_T = Y_0^\top\pi_0 + \int_0^T Y_t^\top \bSig_t H R^{-1} \ud Z_t
\]
where
$$
\frac{\ud Y_t}{\ud t} = (-A + HR^{-1}H^\top \bSig_t) Y_t,\quad Y_T = f
$$ 
Define $\Phi_{t,T}$ to denote the state transition matrix and express the
solution as $Y_t= \Phi_{t,T} f$.  Thus,
\[
S_T = f^\top \underbrace{(\Phi_{0,T}^\top \pi_0 + \int_0^T \Phi_{t,T}^\top\bSig_t H R^{-1} \ud Z_t)}_{=:\bar{X}_T}
\]
Noting that time $T$ is arbitrary, upon differentiating with respect
to $T$, one obtains the Kalman filter
\begin{align*}
\ud \barX_t 
&= A^\top \barX_t \ud t - \bSig_tHR^{-1} (\ud Z_t - H^\top \barX_t\ud t)
\end{align*}
where we have replaced $T$ by $t$.

\subsection{Proof of Thm.~\ref{thm:optimal-solution}}\label{apdx:maximum}
Equation~\eqref{eq:Hamilton_eqns-a}-\eqref{eq:Hamilton_eqns-c} are 
Hamilton's equation for optimal control of a 
BSDE~\cite{peng1993backward}. Explicitly, the partial derivatives are
as follows:
\begin{align*}
{\cal H}_p &= -Ay-Hu-\dv(Hv^\top)+vH^\top\pi\\
{\cal H}_y &= -A^\top p - \pi(Q) y\\
{\cal H}_v &= -\dv^\dagger(p)H + p\pi^\top H - \pi u^\top R - \dv^\dagger(\pi)vR
\end{align*}
Using these formulae, the explicit form of Hamilton's equation is
as follows:
\begin{subequations}\label{eq:Hamilton_eqns_explicit}
	\begin{align}
\ud P_t &= \big(A^\top P_t + \pi_t(Q)Y_t\big)\ud
	t+ \big(\dv^\dagger(P_t) - P_t\pi_t^\top)HR^{-1}\ud \innov_t   \nonumber\\
&\quad + (\pi_tU_t^\top+\dv^\dagger(\pi_t)V_t) \ud \innov_t
	\label{eq:Hamilton_eqns_explicit-a}\\
\ud Y_t &= \big(-AY_t - HU_t - \dv(HV_t^\top) + V_tH^\top \pi_t\big)\ud t + V_t\ud \innov_t \label{eq:Hamilton_eqns_explicit-b}\\
P_0 &= \Sigma_0 Y_0 , \quad Y_T = f \label{eq:Hamilton_eqns_explicit-c}
	\end{align}
\end{subequations}
The optimal control is obtained from the maximum principle:
$$
U_t = \mathop{\mathrm{argmax}}_{\alpha\in\Re^m} \; {\cal
  H}(Y_t,V_t,\alpha,P_t\,;\pi_t)
$$
Since ${\cal H}$ is quadratic in the control input, the explicit
formula~\eqref{eq:opt-cont-soln} is obtained by setting the derivative
to zero:
$$
{\cal H}_u(Y_t,V_t,U_t,P_t\,;\pi_t) = H^\top P_t + RU_t + RV_t^\top \pi_t = 0
$$

\subsection{Proof of Thm.~\ref{thm:Wonham}}\label{apdx:Wonham}

As noted in   Remark~\ref{rm:linear-structure}, the optimal control problem for
the finite state-space case has a linear structure, and thus 
$
P_t = M_t Y_t
$ for some matrix-valued process $\{M_t \in \Re^{d\times d}:0\le t \le T\}$.
The proof is broken into two steps: In step 1, we assume $M_t =
\Sigma_t$ and derive the equation~\eqref{eq:Wonham} of the Wonham
filter. In step 2, we show that this assumption is consistent with the
filter.

\spm{The proof approach is to use Hamilton's equation~\eqref{eq:Hamilton_eqns_explicit-b} for $Y_t$ to derive the equation for $\pi_t$.
\\
This sentence is a bit confusing (proof of what).  Not sure it is needed.}

\newP{Step 1}  
The foregoing implies the following identities  
\begin{equation}
\label{eq:pi_t-y_t}
\begin{aligned}
 \big(P_t^\top H R^{-1} +  &\pi_t^\top V_t\big) \ud Z_t
 \\
 &= \ud(\pi_t^\top Y_t)
             = Y_t^\top \ud \pi_t + \pi_t^\top \ud Y_t + \ud Y_t^\top\ud \pi_t
\end{aligned} 
\end{equation}
The first equality is  \eqref{eq:pi_t_y_t_est}, and the second follows from
   It\^{o}'s product formula.   
 
Hamilton's equation~\eqref{eq:Hamilton_eqns_explicit-b} gives a formula for  $\ud Y_t$,  which when combined with  \eqref{eq:pi_t-y_t} gives
\spm{I thought I had clarified this, and now I see I have made it more obscure.  It does need clarification, but I am out of time}
\begin{align}
Y_t^\top \big(&\ud \pi_t - (A^\top \pi_t\ud t + \Sigma_t^\top HR^{-1}\ud \innov_t)\big)\nonumber \\
&= \big(\pi_t^\top \dv(HV_t^\top) - \pi_t^\top V_t H^\top \pi_t\big)\ud t-(V_t\ud \innov_t)^\top\ud \pi_t \label{eq:dpi_expression}
\end{align}
Upon integrating both sides, one finds that
\[
G_t:=\int_0^t Y_\tau^\top (\ud \pi_\tau - \Sigma_\tau^\top HR^{-1}\ud \innov_\tau)
\]
is a finite variation process.  Therefore, by an application of~\cite[Theorem~4.8]{le2016brownian}, 
$$
\ud \pi_t =g_t\ud t + \Sigma_t HR^{-1}\ud \innov_t
$$
for some yet to be determined process $g$.
Now,
\begin{align*}
(V_t\ud \innov_t)^\top\ud \pi_t &= \tr(V_t^\top \ud \pi_t \ud \innov_t^\top)\\
&=\tr(V_t^\top \Sigma_t H)\ud t\\
&=\big(\tr(\dv^\dagger(\pi_t)HV_t^\top) - \pi_t^\top HV_t^\top \pi_t\big)\ud t\\
&=\big(\dv(HV_t^\top)^\top\pi_t - \pi_t^\top V_t H^\top \pi_t\big)\ud t
\end{align*}
Therefore, the right-hand side of~\eqref{eq:dpi_expression} must be zero:
$$
Y_t^\top \big(\ud \pi_t - (A^\top \pi_t\ud t + \Sigma_t^\top HR^{-1}\ud \innov_t)\big) = 0
$$
This gives the equation of the Wonham filter.

\newP{Step 2} It remains to verify that $P_t = \Sigma_t Y_t$. Using
the equation of the Wonham filter~\eqref{eq:Wonham} and the definition~\eqref{eq:pi_and_I-b} for $\Sigma_t$, it is a direct calculation to see that 
\begin{align}
\ud \Sigma_t
=&\big(A^\top\Sigma_t + \Sigma_t A + \pi_t(Q) -\Sigma_tHR^{-1}H^\top\Sigma_t\big)\ud t \label{eq:DRE-for-NLF}\\
&+\dv^\dagger(\Sigma_t H R^{-1} \ud \innov_t) - \Sigma_t H R^{-1} \ud \innov_t \pi_t^\top- \pi_t\ud \innov_t^\top R^{-1}H^\top\Sigma_t \nonumber
\end{align}
The assertion is shown by establishing
\begin{equation}\label{eq:thm2-assertion}
\ud(\Sigma_t Y_t) = \Sigma_t \ud Y_t + \ud \Sigma_t Y_t + \ud \Sigma_t\ud Y_t = \ud P_t
\end{equation}
and noting $\Sigma_0Y_0 = P_0$.

The calculation showing~\eqref{eq:thm2-assertion} is notationally
cumbersome but straightforward.  It is included in step 3 below.

\newP{Step 3} 
For any two column vectors $a,b\in\Re^d$, $a\cdot b$ denotes the Hadamard
(element-wise) product.  For $a,b\in\Re^d$, it is a straightforward
calculation to see
\begin{align*}
\big(\dv^\dagger(\Sigma_t a)& -\Sigma_t a\pi_t^\top-\pi_ta^\top \Sigma_t\big)b\\
& =\Sigma_t\big(\dv(ab^\top) - ba^\top\pi_t-ab^\top\pi_t\big)
\end{align*}
Multiplying both sides of the matrix-valued
equation~\eqref{eq:thm2-assertion} by $Y_t$, upon using the identity
with $a = HR^{-1}\ud \innov_t$ and $b = {Y}_t$ to simplify the
righthand-side, one obtains 
\begin{align*}
\ud \Sigma_t \; {Y}_t=&\big(A^\top\Sigma_t + \Sigma_t A + \pi_t(Q) -\Sigma_t HR^{-1}H^\top\Sigma_t\big){Y}_t\ud t\\
&+\big(\dv^\dagger(\Sigma_t {Y}_t) - \Sigma_t {Y}_t \pi_t^\top- \pi_t{Y}_t^\top\Sigma_t\big) H R^{-1} \ud \innov_t
\end{align*}
Similarly, multiplying both sides of~\eqref{eq:thm2-assertion} by
${V}_t\ud \innov_t$, using the identity with $a = HR^{-1}\ud \innov_t$
and $b = {V}_t\ud \innov_t$, after applying It\^{o} rules to simplify
the righthand-side, one obtains
\begin{align*}
\ud \Sigma_t {V}_t\ud \innov_t=&\Sigma_t\big(\dv(H{V}_t^\top)-H{V}_t^\top\pi_t - {V}_tH^\top\pi_t\big)\ud t
\end{align*}
Therefore,
\begin{align*}
\ud (\Sigma_t Y_t) & = \Sigma_t \big(-AY_t +HR^{-1}H^\top \Sigma_tY_t\big)\ud t\\
& \quad + \Sigma_t \big(HV_t^\top\pi_t + V_tH^\top\pi_t - \dv(HV_t^\top)\big)\ud t\\
& \quad + \Sigma_t V_t\ud \innov_t + \ud \Sigma_t Y_t +\ud \Sigma_t V_t\ud \innov_t\\
&= \big(A^\top \Sigma_t Y_t + \pi_t(Q)Y_t\big)\ud t + \Sigma_t V_t\ud \innov_t \\
& \quad + \big(\dv^\dagger(\Sigma_t Y_t)-\Sigma_t Y_t \pi_t^\top - \pi_t Y_t^\top \Sigma_t\big)HR^{-1}\ud \innov_t
\end{align*}
The right-hand side is identical to the right-hand side of Hamilton's equation~\eqref{eq:Hamilton_eqns_explicit-a} for $P_t$.  

\subsection{Proof of Thm.~\ref{thm:optimal-solution-euc}}\label{apdx:L2-derivatives}

Equation~\eqref{eq:Hamilton_eqns_euc} are Hamilton's equation.
The Gateaux differentials of ${\cal H}$ are  
\begin{align*}
{\cal H}_p &= -\clA y - h^\top (u + v) + v^\top \pi(h)\\
{\cal H}_y &= -\clA^\dagger p + \operatorname{div} \big(\pi\sigma\sigma^\top \frac{\partial y}{\partial x} \big)\\
{\cal H}_v &= -ph + p\langle h,\pi\rangle - R(u+v)\pi\\
{\cal H}_u &= \langle -h, p\rangle - Ru - \langle Rv,\pi\rangle
\end{align*}
Therefore the explicit formulas for Hamilton's equations are
\begin{align*}
\ud P_t(x) &=\Big((\clA^\dagger P_t)(x) - \operatorname{div} \big(\pi_t \sigma \sigma^\top \frac{\partial Y_t}{\partial x}\big)(x) \Big)\ud t 
\\
&\ \
	+ \big(P_t(x) (h(x) -\langle\pi_t , h \rangle)^\top R^{-1} + \pi_t(x)(U_t+V_t(x))^\top\big)\ud \innov_t
\\
\ud Y_t(x) &= \big(-(\clA Y_t)(x) - h^\top(x) (U_t + V_t(x)) + V_t^\top(x) \langle \pi_t, h\rangle \big)\ud t 
\\
&
\ \ 
	+ V_t(x) \ud \innov_t
\end{align*}
The optimal control is obtained from the maximum principle:
$$
U_t = \mathop{\mathrm{argmax}}_{\alpha\in\Re^m} \; {\cal H}(Y_t, V_t, \alpha, P_t \,; \pi_t)
$$
The explicit formula~\eqref{eq:opt-cont-soln-euc} for the optimal
control is obtained by setting the derivative to zero:
$$
{\cal H}_u(Y_t,V_t,U_t,P_t\,;\pi_t) = \langle -h, P_t\rangle - RU_t - \langle RV_t,\pi_t\rangle = 0
$$

\medskip
\subsection{Proof of Thm.~\ref{thm:Kushner}}\label{apdx:Kushner}

The proof is identical to the proof of Theorem~\ref{thm:Wonham} for
the finite state-space case.  In step 1, we derive the equation of the
filter by first assuming the following linear relationship between $P_t$ and $Y_t$:
\begin{equation}\label{eq:P_t-conjecture}
P_t(x) = \pi_t(x)\big(Y_t(x)-\langle \pi_t,Y_t\rangle\big)
\end{equation}
In step 2, we verify this relationship is consistent with the filter equation.

\newP{Step 1} Suppose~\eqref{eq:P_t-conjecture} is true. The
differential form of~\eqref{eq:pi_t-y_t-euc} is given by:
\begin{equation}\label{eq:ytpitdiffeuc}
\ud\big(\langle \pi_t, Y_t\rangle \big) = \big(R^{-1}\langle P_t,h\rangle + \langle \pi_t, V_t\rangle \big)^\top \ud Z_t
\end{equation}
As in the finite state-space case, the proof approach is to use Hamilton's equation for $Y_t$ to derive the equation for $\pi_t$.

Using the It\^{o} product formula,
\begin{align*}
\ud\big(\langle \pi_t, Y_t\rangle \big)
&= \langle \ud \pi_t, Y_t\rangle + \langle \pi_t, \ud Y_t\rangle + \langle \ud \pi_t, \ud Y_t\rangle
\end{align*}
Use Hamilton's equation~\eqref{eq:Hamilton_eqns_euc-b} to evaluate $\ud Y_t$, and equate
the resulting expression to the right-hand side
of~\eqref{eq:ytpitdiffeuc} to obtain:
\begin{align*}
\langle Y_t, \ud \pi_t & - (\clA^\dagger\pi_t \ud t + \pi_t(h-\langle \pi_t,h\rangle )^\top R^{-1} \ud \innov_t)\rangle \\
&= -\langle V_t^\top\ud \innov_t, \ud \pi_t \rangle +\langle V_t^\top, \pi_t(h - \langle \pi_t,h\rangle)\rangle \ud t
\end{align*}
This is the Euclidean counterpart of~\eqref{eq:dpi_expression} in the
proof of Theorem~\ref{thm:Wonham} in the finite state-space case.  The
derivation of the filter is now identical.  

\newP{Step 2} The verification of~\eqref{eq:P_t-conjecture} follows
along the same lines as the finite state-space case.  It is omitted
here.

\medskip
\subsection{Proof of Thm.~\ref{thm:martingale-DP}}\label{apdx:martingale}

The proof is given for the finite state-space case. In the finite state-space case,
$$
{\cal V}(y\,;\pi_t) = \half y^\top \Sigma_t y
$$ 
Upon using~\eqref{eq:pi_and_I-b} for $\Sigma_t$ and~\eqref{eq:Wonham} for the filter, it is a direct application of the It\^{o} product formula that:
\begin{align*}
\ud (Y_t^\top\Sigma_t Y_t)
&=\big(Y_t^\top \pi_t(Q)Y_t - Y_t^\top \Sigma_tHR^{-1}H^\top \Sigma_t Y_t+\tr(V_t^\top\Sigma_tV_tR) \big)\ud t\\
&\big(-2U_t^\top H^\top \Sigma_t Y_t-2Y_t^\top \Sigma_t HV^\top \pi_t\big)\ud t +2Y_t^\top \Sigma_t V_t\ud \innov_t\\
&+ Y_t^\top(\dv^\dagger(\Sigma_t Y_t)-\Sigma_t Y_t \pi_t^\top - \pi_t Y_t^\top \Sigma_t)HR^{-1}\big)\ud \innov_t
\end{align*}
Using this formula,
\begin{align*}
\ud &{\cal M}_t^U = \half \ud (Y_t^\top\Sigma_t Y_t) - {\cal L}(Y_t,V_t,U_t\,;\pi_t)\ud t\\
&=-\half(U_t + R^{-1} H^\top \Sigma_t Y_t + V_t^\top \pi_t)^\top R(U_t + R^{-1} H^\top \Sigma_t Y_t + V_t^\top \pi_t)^\top \ud t\\
&\quad +Y_t^\top \Sigma_t V_t\ud \innov_t + \half Y_t^\top(\dv^\dagger(\Sigma_t Y_t)-\Sigma_t Y_t \pi_t^\top - \pi_t Y_t^\top \Sigma_t)HR^{-1}\ud \innov_t
\end{align*}
Therefore, ${\cal M}^U$ is always a super-martingale with respect to $\clZ$, and is a martingale if and only if 
$$
U_t = -R^{-1}H^\top \Sigma_t Y_t - V_t^\top \pi_t
$$
for almost every $t \in [0,T]$.
Consequently, 
\[
\E\big({\cal M}_T^U\big) \leq \E\big({\cal M}_0^U\big)
\]
Adding $\E\big(\int_0^T {\cal L}(Y_s,V_s,U_s\,;\pi_s)\ud s\big)$ on both sides yields the optimality result.

\end{document}